\documentclass[a4paper]{article}
\pdfoutput=1

\usepackage{IEK10} 
\usepackage{natbib}

\makeatother
\usepackage{amssymb}
\usepackage{xcolor}
\usepackage{mathtools, cuted}
\usepackage{tabularx}
\usepackage{eurosym}
\usepackage[utf8]{inputenc}

\newcommand{\vect}[1]{\mathbf{#1}}

\renewcommand{\min}[1][]{
	\ifthenelse{\isempty{#1}}{\operatorname{min}}{\ensuremath{\underset{#1}{\text{min}\,}}}
}
\newcommand{\Should}{\ensuremath{\overset{!}{=}}}

\def\IEK10{
  Forschungszentrum Jülich GmbH,
  Institute of Energy and Climate Research,
  Energy Systems Engineering (IEK-10),
  Jülich 52425,
  Germany
}
\def\RWTH{
  RWTH Aachen University
  Aachen 52062,
  Germany
}
\def\ETH{
  ETH Zürich,
  Energy \& Process Systems Engineering,
  Zürich 8092,
  Switzerland
}

\newcommand{\mytitle}{Dynamic Ramping for Demand Response of Processes and Energy Systems based on Exact Linearization}

\newcommand{\affil}{
  \begin{itemize}[leftmargin=3mm, itemsep=0mm]
    \item[$^a$]\IEK10
    \item[$^b$]\RWTH
    \item[$^c$]\ETH
  \end{itemize}
}

\def\firstAuthor{Florian Joseph Baader}
\newcommand{\myauthor}{\firstAuthor$^{a,b,c}$, Philipp Althaus$^{a,b}$, André Bardow$^{a,c}$, Manuel Dahmen$^{a,*}$}
\newcommand{\correspondingauthor}{M. Dahmen, \IEK10 \newline E-mail: m.dahmen@fz-juelich.de}
\author{\myauthor}

\usepackage[
  colorlinks,
  linkcolor=blue,
  citecolor=blue,
  urlcolor=blue,
  pdftitle={\mytitle},
  pdfauthor={\firstAuthor}
]{hyperref}
\usepackage[capitalise, nameinlink]{cleveref}
\crefname{table}{Tab.}{Tab.}

\newcommand{\setpgfexternalcounter}[1]{
  \makeatletter%
  \pgfkeysgetvalue{/tikz/external/figure name}\myexternalname
  \expandafter\gdef\csname c@tikzext@no@\myexternalname\endcsname{#1}%
  \makeatother
}

\begin{document}

  \thispagestyle{firststyle}

  \begin{center}
    \begin{large}
      \textbf{\mytitle}
    \end{large} \\
    \myauthor
  \end{center}

  \vspace{0.5cm}

  \begin{footnotesize}
    \affil
  \end{footnotesize}

  \vspace{0.5cm}

  \begin{abstract}
The increasing share of volatile renewable electricity production motivates demand response.
Substantial potential for demand response is offered by flexible processes and their local multi-energy supply systems. 
Simultaneous optimization of their schedules can exploit the demand response potential, but leads to numerically challenging problems for nonlinear dynamic processes.
In this paper, we propose to capture process dynamics using dynamic ramping constraints. 
In contrast to traditional static ramping constraints, dynamic ramping constraints are a function of the process state and can capture high-order dynamics. 
We derive dynamic ramping constraints rigorously for the case of single-input single-output processes that are exactly input-state linearizable.
The resulting scheduling problem can be efficiently solved as a mixed-integer linear program. 
In a case study, we study two flexible reactors and a multi-energy system. 
The proper representation of process dynamics by dynamic ramping allows for faster transitions compared to static ramping constraints and thus higher economic benefits of demand response.
The proposed dynamic ramping approach is sufficiently fast for application in online optimization.
\end{abstract}

\vspace{0.5cm}


\vspace{0.75cm}

\begin{flushleft}
  \leavevmode{\parindent=15mm\indent}
  \textbf{Highlights:}
  \begin{itemize}[leftmargin=20mm]
    \item Dynamic ramping constraints for scheduling of processes and energy systems 
    \item Allow for high-order dynamics, and non-constant ramp limits
    \item Derived rigorously from exact input-state linearization
    \item Mixed-integer linear optimization fast enough for online optimization
    \item Cost savings close to nonlinear optimization
  \end{itemize}
\end{flushleft}
\vspace*{5mm}

  \newpage

\section{Introduction}
Many countries are transforming  their national energy systems towards renewable energies. 
This  transformation requires more renewable electricity generation by wind and sun.
The inherent volatility of renewable electricity generation causes temporal imbalances of demand and supply. 
These imbalances can be reduced if consumers shift their demand in time.
Ideally, both electricity grid and consumers benefit from this demand response (DR) \citep{Zhang.2016}.
To incentivize demand shifting for consumers, electricity is traded with time-varying prices at day-ahead and intra-day markets.
At these markets prices react to the current demand and supply.  

\begin{figure}[h]
\centering
\includegraphics{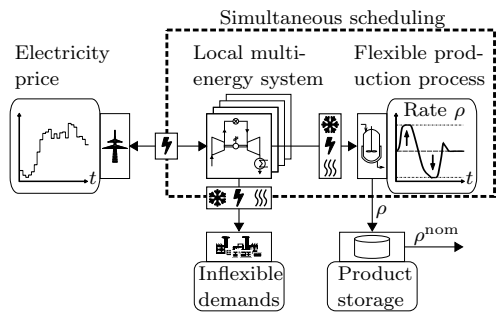}
  \caption{Simultaneous scheduling of a flexible production process and its local multi-energy system reacts to variable electricity prices using the demand response potential of the process by modulating the production rate $\rho$ and thus process energy consumption. }
   \label{fig:system}
\end{figure}

Demand response is especially promising for the chemical industry. Many energy-intensive production processes offer demand response potential, i.e., these processes can (i) adjust their production rate and thus energy demand, and (ii) store their product for later use \citep{Mitsos.2018,Merkert.2015}. 
However, chemical processes consume different forms of energy – not only electricity,  but also cooling, or heating.
Therefore, processes are typically supplied by a local multi-energy system  which consumes primary energy sources and exchanges electricity with the grid, either by buying electricity from the grid or by selling electricity from on-site production \citep{Voll.2013}. 
Consequently, DR needs to consider both the process and its local energy system in a simultaneous scheduling \citep{MujtabaH.Agha.2010,Leenders.2019b,baader2022simultaneous} (Figure \ref{fig:system}) that determines operational set points for a time horizon in the order of one day \citep{Baldea.2014}.

The desired simultaneous scheduling leads to computationally intensive optimization problems because processes and energy systems introduce two different challenges that are hard to solve simultaneously:  First, processes often introduce scheduling-relevant nonlinear dynamics \citep{Mitsos.2018,Baldea.2014,Daoutidis.2018,Caspari.2019_2,Otashu.2019}. 
Scheduling-relevant means the time that the process needs to change between two steady states is in the same order of magnitude as electricity price time steps, e.g., 1 hour.
Second, energy system models often introduce discrete decision variables resulting from often redundant units with minimum part-load constraints \citep{Voll.2013}. 
Thus, the resulting problems are mixed-integer nonlinear dynamic optimization problems which are notoriously difficult to solve. 
However, scheduling optimization must be performed online in order to timely provide operational set-points for the underlying control. 
The maximum allowed optimization runtime is typically between 5 and 20 minutes \citep{IiroHarjunkoski.2014}. 
To achieve such runtimes, the optimization problem needs to be reformulated to a simpler problem.

Online-applicable optimization runtimes of local energy systems alone can often be achieved by mixed-integer linear programming (MILP) formulations \citep{MichaelJ.Risbeck.2017,Mitra.2013,Carrion.2006,SusanneSass.2020}. 
The MILP problem can also integrate the production process if the process has negligible dynamics such that quasi-steady-state assumptions can be used on the scheduling time scale. 
Scheduling can then calculate a feasible trajectory of operational set-points. 
Based on this trajectory, the real process inputs are determined by the underlying control \citep{Baldea.2014}. 
Accordingly, for negligible process dynamics, no process model is needed but only the process energy demands must be described as piece-wise affine function of the production rate \citep{Schafer.2020,Bree.2019}. 

For processes with relevant dynamics, scheduling optimization still does not necessarily need to consider the full-order process model but rather requires a set of constraints that determine how fast the process can change its production rate. 
Traditional first-order ramping constraints bound the first derivative of the production rate using constant limits \citep{Carrion.2006} but have two shortcomings for simultaneous scheduling that we address in this paper: 
First, constant ramping limits cause either unnecessarily conservative or infeasible schedules if the achievable rate of change varies with the process state \citep{T.Li.2007,C.M.CorreaPosada.2017}, which is typical for nonlinear chemical processes. 
Second, first-order dynamics may not be applicable.
If, for example, a chemical reactor has a temperature-dependent production rate and is cooled through a reactor jacket with a significant thermal inertia, changing the production rate through a change of reactor temperature would lead to at least second-order dynamics. 

To overcome both shortcomings, we propose high-order dynamic ramping constraints for simultaneous DR scheduling of processes and energy systems. 
For this purpose, we present a method to rigorously derive such dynamic ramping constraints for the case of exact input-state linearizable single-input single-output (SISO) processes based on the full-order nonlinear process model. 
Our high-order dynamic ramping constraints are based on a chain of differential equations. 
The highest considered time-derivative is the ramping degree of freedom that is limited by dynamic limits as function of the process state. 
With time discretization by collocation  \citep{Biegler.2010}, these high-order dynamic ramping constraints can be converted to linear algebraic constraints and thus allow for an MILP formulation.
Thus, dynamic ramping constraints can be readily integrated into typical MILP-based energy system optimization models \citep{SusanneSass.2020}.

The remaining paper is structured as follows: In Section \ref{sec:prolem_formulations}, we introduce the original nonlinear optimization problem for simultaneous scheduling and the reformulated linear optimization problem with dynamic ramping constraints. 
In Section \ref{sec:method}, we present a rigorous derivation of  dynamic ramping constraints from exact input-state linearization for SISO processes. 
In Section \ref{sec:case_study}, a case study featuring two continuous stirred tank reactors (CSTRs) and a multi-energy system is investigated. 
In Section \ref{sec:discussion}, we discuss possible extensions for cases that are not covered by our current assumptions.
Section \ref{sec:conclusion} concludes the work.

\section{Simultaneous dynamic scheduling of process and energy system}
\label{sec:prolem_formulations}
In Section \ref{sec:P1}, we present the original simultaneous dynamic scheduling problem (P1) which is a nonlinear mixed-integer dynamic optimization (MIDO) problem.
In Section \ref{sec:P2}, we introduce the proposed MILP problem formulation (P2) based on high-order dynamic ramping constraints.

\subsection {Nonlinear mixed-integer dynamic scheduling}
\label{sec:P1}
In the original scheduling optimization problem (P1), all decision variables
        $\boldsymbol{\chi} = (\vect{x}^T, \vect{u}^T, (\vect{Q}_{\text{dem}}^{\text{process}})^T, \rho, S,$ $\Phi_{\text{energy}},   (\vect{Q}^{\text{in}})^T, (\vect{Q}^{\text{out}})^T, (\Delta\vect{P})^T, \vect{z}_{\text{on}}^T )^T $,
which are further explained in the following, are functions of time $t$ although not stated explicitly to improve readability. We use the notation $\dot{\chi}$ to indicate the first time derivative of a variable $\chi$, and  $\chi^{(k)}$ to indicate the $k$-th time derivative.
The problem reads:
    \begin{align}
        \label{eq:P1_objective} \tag{P1a}
        &\hspace{2cm}  \underset{\boldsymbol{\chi} \in \left[\boldsymbol{\chi}^l, \boldsymbol{\chi}^u\right]}{\text{min}}~~  \Phi_{\text{energy}}(t_f) \\
        \text{s.t. }~& \text{Process model: } \nonumber
        \\ \label{eq:P1_process_model} \tag{P1b}
        & \dot{\vect{x}}=\vect{f}_1(\vect{x}) + \vect{f}_2(\vect{x})\vect{u} ~~~\forall t \in \left[t_0, t_f\right]\\
        & \text{Product quality: } \nonumber
        \\
        \label{eq:P1_product_qual} \tag{P1c} & 0 \geq \vect{g}(\vect{x}, \rho) ~~~\forall t \in \left[t_0, t_f\right]
        \\
        \nonumber
        & \text{Process energy demand: ~}
        \\ \label{eq:P1_process_energy_dem} \tag{P1d}& Q_{\text{dem},e}^{\text{process}} = h_e(\vect{x},\vect{u}) ~\forall e \in \mathbb{E}, ~\forall t \in \left[t_0, t_f\right]
        \\
        & \text{Product storage: } \nonumber
        \\ \label{eq:P1_storage} \tag{P1e}&  \dot{S} = \rho - \rho^{\text{nom}} ~~~\forall t \in \left[t_0, t_f\right]
        \\  \nonumber
        & \text{Energy costs: }
        \\ \label{eq:P1_energy_costs} \tag{P1f}& \dot{\Phi}_{\text{energy}} = \sum_{e \in \mathbb{E}} p_e\left( \sum_{i \in \mathbb{C}_e^{\text{cons}}}  Q_{i}^{\text{in}} + \Delta P_e \right) ~~~\forall t \in \left[t_0, t_f\right]
        \\
        \nonumber
        & \text{Energy conversion: ~}
        \\ \label{eq:P1_eff} \tag{P1g}& Q_{i}^{\text{out}} = \eta_i \left( \frac{Q_{i}^{\text{out}}}{Q_{i}^{\text{max}}} \right)Q_{i}^{\text{in}}  ~~~\forall i \in \mathbb{C},~\forall t \in \left[t_0, t_f\right]
        \\ 
        \nonumber
        & \text{Minimum part-load: }
        \\ \label{eq:P1_part_load} \tag{P1h}
        & z_i^{\text{on}} Q_i^{\text{min}} \leq Q_i^{\text{out}} \leq  z_i^{\text{on}} Q_i^{\text{max}} ,~~~\forall i \in \mathbb{C},~\forall t \in \left[t_0, t_f\right]
        \\ \nonumber &\text{with} ~ z_i^{\text{on}} \in \{0,1\}
        \\ 
        \nonumber
        & \text{Energy balance: }
        \\ \label{eq:P1_energy balance} \tag{P1i}& 
        Q_{\text{dem},e}^{\text{process}} + Q_{\text{dem},e}^{\text{inflexible}} = \sum_{i \in \mathbb{C}_e^{\text{sup}}} Q_{i}^{\text{out}} +\Delta P_e
        ~~~\forall e \in \mathbb{E},  ~\forall t \in \left[t_0, t_f\right] 
    \end{align}

The objective (\ref{eq:P1_objective}) is to minimize the cumulative energy costs $\Phi_{\text{energy}}$ at final time $t_f$. 
All decision variables are subject to upper and lower bounds, $\boldsymbol{\chi}^u, \boldsymbol{\chi}^l$, respectively.
The process model (\ref{eq:P1_process_model})  of chemical production processes can usually be expressed in input-affine control form \citep{Corriou.2018,baldea_daoutidis_2012}, i.e., the time derivative of process states $\dot{\vect{x}}$ is given by two nonlinear functions $\vect{f}_1(\vect{x})$, $\vect{f}_2(\vect{x})$ and the process degrees of freedom $\vect{u}$. The process model is valid for all time points between initial time $t_0$ and final time $t_f$.
To maintain product quality (\ref{eq:P1_product_qual}), we assume constraints $\vect{g}$ on process states $\vect{x}$ and the production rate $\rho$.
The process energy demand  $Q_{\text{dem},e}^{\text{process}}$ (\ref{eq:P1_process_energy_dem}) for an energy form $e$ in the set of energy forms $\mathbb{E}$ is a function $h_e(\vect{x},\vect{u})$ of process states $\vect{x}$ and degrees of freedom $\vect{u}$.
For product storage (\ref{eq:P1_storage}), we assume a buffer storage with filling level $S$.
As all decision variables, the filling level $S$ is subject to upper and lower bounds.
Moreover, the final storage filling level $S(t_f)$ might be constrained to be greater than or equal to the initial filling level $S(t_0)$ to avoid depletion of the storage \citep{Schafer.2020}.
The storage unit is filled by the production rate $\rho$ of the process and emptied with a constant nominal product demand rate $\rho^{\text{nom}}$ \citep{Caspari.2019_2,Schafer.2020,Pattison.2016}.
The instantaneous energy costs (\ref{eq:P1_energy_costs})  are the sum over specific price times consumption for all energy forms in $\mathbb{E}$ given by the 
input power of energy system components $Q_{i}^{\text{in}}$,
the set of energy system components that consume energy $e$, $\mathbb{C}_e^{\text{cons}}$, the energy prices $p_e$,  and the power exchanged with the grid $\Delta P_e$.
For the energy conversion (\ref{eq:P1_eff}), the output power $Q_{i}^{\text{out}}$ of each component $i$ in the set of components $\mathbb{C}$ is the product of the input power $Q_{i}^{\text{in}}$ and the efficiency $\eta_i$. 
The efficiency $\eta_i$ itself is a function of the part-load fraction, i.e., the output power $Q_{i}^{\text{out}}$ divided by the maximum output power $Q_i^{\text{max}}$ \citep{SusanneSass.2020}.
The minimum part-load constraints (\ref{eq:P1_part_load}) for energy system components $i$ require a binary variable $z_i^{\text{on}}$ which ensures that if the component is on the output power is between maximum and minimum value $Q_i^{\text{max}}$ and $Q_i^{\text{min}}$, respectively \citep{Voll.2013}. 
The energy balance (\ref{eq:P1_energy balance}) states that for every energy form $e$ the demands of the flexible production process $Q_{\text{dem},e}^{\text{process}}$ and other inflexible processes $Q_{\text{dem},e}^{\text{inflexible}}$ (cf. Figure \ref{fig:system}) must be met by the set of energy system components that supply $e$, $\mathbb{C}_e^{\text{sup}}$. 
Additionally, power $\Delta P_e$ can be exchanged with the electricity grid.
Finally, the initial values $\vect{x}_0$ and $S_0$, provide the initial conditions $\vect{x}(t_0) = \vect{x}_0$, $S(t_0) = S_0$, and $\Phi_{\text{energy}}(t_0) = 0$. 
In the original optimization problem (P1), we find differential Equations  (\ref{eq:P1_process_model}, \ref{eq:P1_storage}, \ref{eq:P1_energy_costs}), nonlinear Equations (\ref{eq:P1_process_model}, \ref{eq:P1_product_qual}, \ref{eq:P1_process_energy_dem}, \ref{eq:P1_eff}), and one binary variable ($z_i^{\text{on}}$) per energy system component and timestep in Equation (\ref{eq:P1_part_load}).
To solve the optimization problem in online-applicable runtime, (P1) needs to be simplified.

\subsection{Linear mixed-integer dynamic scheduling with ramping constraints}
\label{sec:P2}
The energy system part of problem (P1) can be reformulated as an MILP if nonlinear efficiency curves (\ref{eq:P1_eff}) are approximated by piece-wise affine functions \citep{SusanneSass.2020}.
To integrate the process into a MILP formulation, the nonlinear process model (\ref{eq:P1_process_model} - \ref{eq:P1_process_energy_dem}) must be replaced by simpler linear equations.
As discussed in the introduction, these simpler equations must capture the information on how fast the production rate $\rho$ of the process can be changed and must provide an sufficiently accurate approximation of the process energy demand.
These requirements can be fulfilled by a combination of piece-wise affine (PWA) ramping constraints plus a PWA process energy demand model. 
This combination replaces (\ref{eq:P1_process_model}) - (\ref{eq:P1_process_energy_dem}).
Consequently, the vectors $\vect{x}$ and $\vect{u}$, which only occur in the removed Equations (\ref{eq:P1_process_model}) - (\ref{eq:P1_process_energy_dem}), are removed from the optimization variables $\boldsymbol{\chi}$ and the problem (P2) is given by:
    \begin{align}
        \label{eq:P2_objective} \tag{P2a}
        &\hspace{0.5cm}  \underset{\boldsymbol{\chi} \in \left[\boldsymbol{\chi}^l, \boldsymbol{\chi}^u\right]}{\text{min}}~~  \Phi_{\text{energy}}(t_f) \\
        \text{s.t.}~~\label{eq:P2_dynamic_ramping} \tag{P2b}
        & \text{PWA ramping constraints}~ \forall t \in \left[t_0, t_f\right]
        \\
        \label{eq:P2_process_energy_dem} \tag{P2c}
        & \text{PWA process energy demand model} ~\forall t \in \left[t_0, t_f\right]
        \\ 
        \nonumber & \text{Equations (\ref{eq:P1_storage}) - (\ref{eq:P1_energy balance})}
    \end{align}
For (\ref{eq:P2_dynamic_ramping}), traditional first-order ramping constraints with static limits could be used \citep{Carrion.2006}. Here, static means that the bounds do not depend on the process state.
Such first-order static ramping constraints (SRC) use the first time derivative of the production rate as ramping degree of freedom $\nu$ \citep{Carrion.2006}: 
    \begin{gather}
        \label{eq:ramp_constr} \tag{SRCa}
        \dot{\rho} = \nu,~~~\text{with} \\ \tag{SRCb} \nu^{\text{min}}\leq\nu\leq\nu^{\text{max}}
    \end{gather}
In SRC, the degree of freedom $\nu$  is bounded by constant limits $\nu^{\text{min}}$, $\nu^{\text{max}}$.
By applying time discretization, the static ramping constraints can be converted to linear algebraic constraints.
However, the restriction to constant limits $\nu^{\text{min}}$, $\nu^{\text{max}}$ often enforces a conservative parameterization for nonlinear chemical processes where the achievable rate of change is in general non-constant.
Moreover, first-order dynamics might not be applicable, e.g., if it is necessary to first overcome a temperature inertia before changing the production rate.

To overcome both shortcomings of static ramping constraints, we propose high-order dynamic ramping constraints (DRC).
Dynamic ramping constraints define the $\delta$-th derivative of the production rate as ramping degree of freedom $\nu$ such that the highest time derivative acts as the free variable in scheduling optimization.
Accordingly, the highest derivative can be chosen in every time step and the other derivatives result from time integration.
The ramping degree of freedom $\nu$ is bounded with multivariate functions that depend on the production rate and its time derivatives:
\begin{gather}
     \rho^{(\delta)} = \nu \tag{DRCa} \label{eq:dyn_ramping_a},~~~\text{with}\\  \nu^{\text{min}}\left(\rho,\dot{\rho},...,\rho^{(\delta-1)}\right) \leq \nu \leq \nu^{\text{max}}\left(\rho,\dot{\rho},...,\rho^{(\delta-1)} \right)  \tag{DRCb} \label{eq:dyn_ramping_b}
\end{gather}

The process energy demand (\ref{eq:P2_process_energy_dem}) is often modeled as a piece-wise affine function of the production rate \citep{Schafer.2020}.
We use a multivariate function of the production rate and its time derivatives because during transient operation the energy demand can depend on the speed of the transition:
\begin{align}
    \label{eq:energy_demand_model}
    Q_{\text{dem},e}^{\text{process}} = h_e\left(\rho,\dot{\rho},...,\rho^{(\delta-1)},\nu\right) ~\forall e \in \mathbb{E}
\end{align}
To discretize problem (P2), we use orthogonal collocation on finite elements as an accurate discretization method requiring relatively few discretization points \citep{Biegler.2010}.
As the differential equations introduced by the dynamic ramping constraint (\ref{eq:dyn_ramping_a}), the storage model (\ref{eq:P1_storage}), and the energy costs (\ref{eq:P1_energy_costs}) are all linear, a discretization with collocation in discrete time leads to linear constraints.
If additionally  linear or piece-wise affine approximations are chosen for the limits $\nu^{\text{min}}$ and $\nu^{\text{max}}$ (both in Equation (\ref{eq:dyn_ramping_b})), the energy demand $h_e$ (in Equation (\ref{eq:energy_demand_model})), and the nonlinear  efficiencies $\eta_c$ (in Equation (\ref{eq:P1_eff})), the entire problem P2 can be formulated as MILP.

The problem formulation (P2) allows to model the flexibility of chemical production processes more accurately compared to first-order static ramping constraints while reducing the computational complexity compared to the original nonlinear MIDO problem (P1).
The dynamic ramping constraints are parameterized by the order $\delta$ and the limits $\nu^{\text{min}}$, $\nu^{\text{max}}$ as functions of the process state.
This parameterization could in principle be done based on intuition or based on a suitable heuristic. 
In the following chapter, we show that dynamic ramping constraints can be derived rigorously for the special case of exact input-state linearizable single-input single-output (SISO) processes.

\section{Deriving rigorous dynamic ramping constraints}
\label{sec:method}
In this section, we use the concept of exact linearization from nonlinear control \citep{Corriou.2018} to rigorously derive dynamic ramping constraints.
The derivation is restricted to single-input single-output (SISO) processes that are exact input-state linearizable \citep{Corriou.2018}. 
Before the derivation in Section~\ref{subsec:derive_dyn_ramping}, we state our assumptions in Section~\ref{subsec:assumptions}.
In Section \ref{subsec:PWL_limits}, we discuss piece-wise affine approximations of the true nonlinear ramping limits and the trade-off between conservativeness and computational burden.

\subsection{Assumptions}
\label{subsec:assumptions}

\begin{enumerate}
    \item The process can be described in input-affine control normal form \citep{Corriou.2018,baldea_daoutidis_2012}.
    \item The process has exactly two degrees of freedom: the control input $u$ and the variable production rate $\rho$, i.e., the process model in input-affine control normal form is given by
    \begin{align}
        \label{eq:process_model}
        \dot{\vect{x}} = \vect{f}_1 (\vect{x}) + \left( \vect{f}_{2,1}(\vect{x}), \vect{f}_{2,2}(\vect{x}) \right) \left( \begin{array}{c} u \\ \rho \end{array} \right)~,
    \end{align}
    with state vector $\vect{x} \in \mathbb{R}^n$, and nonlinear functions $\vect{f}_1 (\vect{x})$, $\vect{f}_{2,1} (\vect{x})$, $\vect{f}_{2,2} (\vect{x})$ (compare to \ref{eq:P1_process_model}).
    \item The trajectory of the production rate $\rho$ is determined on the scheduling level. Consequently, on the control level, the input $u$ is the only degree of freedom.
    \item There exists an process output $y$ that is relevant for product quality and that should be maintained constant at its nominal value $y^{\text{nom}}$. 
    This process output can be expressed as a function of the states, $h(\vect{x})$, i.e., the quality constraint \ref{eq:P1_product_qual} simplifies from $0\geq \vect{g}(\vect{x},\rho)$ to $h(\vect{x}) = y^{\text{nom}}$.
    \item The output $y=h(\vect{x})$ can be controlled by exact input-state linearization using the input $u$ \citep{Corriou.2018}. In simple terms, this input-state linearizability is given if the number of inertias between input $u$ and output $y$ is equal to the number of process states $n$ and the state vector $\vect{x}$ can be given as function of the output $y$ and its first $(n-1)$ time derivatives \citep{Corriou.2018}.  
    \item The input $u$ is bounded by a minimum and a maximum value $u^{\text{min}}$, $u^{\text{max}}$, respectively. Both values $u^{\text{min}}$ and $u^{\text{max}}$ are assumed to be constant.
\end{enumerate}

\subsection{Deriving nonlinear ramping limits based on exact linearization}
\label{subsec:derive_dyn_ramping}
In this section, we use exact linearization \citep{Corriou.2018,Slotine.1991,Isidori.1995} to calculate with which dynamics of the production rate $\rho$ control can still hold the output $y$ at the nominal value $y^{nom}$.
From this analysis, we determine dynamic ramping constraints that scheduling optimization has to obey when choosing the trajectory of the production rate $\rho$. 

The output $y$ can be held constant if there always exists a value of the input $u$ that sets the derivatives of $y$ to zero and at the same time is within the bounds $u^{\text{min}}$, $u^{\text{max}}$.
Note that our analysis does not depend on the particular type of control employed because we only analyze if a suitable input $u$ can be chosen on the control level in principle.
Consequently, the derived ramping constraints are a property of the process in combination with the limits on the input $u$. 

On the control level, the process (Equation \ref{eq:process_model}) is a SISO process with the input $u$ and the disturbance $\rho$ that is known in advance.
Note that we consider the production rate $\rho$ to be a disturbance on the control level as it follows the trajectory determined on the scheduling level and thereby induces transient operation. 
To ease notation of the disturbance $\rho$ and its derivatives, we introduce the ramping state vector 
    \begin{align}
        \boldsymbol{\varphi} = \left( \begin{array}{c} \rho \\ \dot{\rho} \\\vdots \\\rho^{(\delta-2)}\\\rho^{(\delta-1)}\end{array} \right) \text{  and its time derivative } \dot{\boldsymbol{\varphi}} = \left( \begin{array}{c} \dot{\rho} \\ \rho^{(2)} \\\vdots \\\rho^{(\delta-1)}\\\nu \end{array} \right) .
    \end{align}

To compensate the disturbance introduced by varying the production rate $\rho$, control manipulates the input $u$.
Thereby, for a process fulfilling our assumptions, control acts on the $r$-th derivative of the output $y$.
The number $r$ is the relative degree defined as the number of times the output $y=h(\vect{x})$ has to be differentiated with respect to time until the input $u$ appears explicitly.
As the process is input-state linearizable (assumption 4), the relative degree $r$ is equal to the number of states $n = \text{dim}(\vect{x})$, i.e., $r=n$ \citep{Corriou.2018}.
Performing $(n-1)$ time differentiations of $y$ gives the first $(n-1)$ derivatives of $y$ as nonlinear functions $\alpha_k(\vect{x},\boldsymbol{\varphi})$ of process states $\vect{x}$ and ramping state vector $\boldsymbol{\varphi}$ with $k=0,...,n-1$.
Our assumption $r=n$ implies that the first $(n-1)$ derivatives of the output $y$ do not depend on $u$, i.e., the term $\frac{\partial \alpha_k(\vect{x},\boldsymbol{\varphi})}{\partial\vect{x}}\vect{f}_2 (\vect{x})$ is equal to zero for all $0<k<n$ (see, e.g., Equation~\ref{eq:y_1}).

By calculating the total time derivative, the derivatives of the output $y$  read:
    \begin{align}
        \tag{4.1}
        y &= h(\vect{x})\coloneqq\alpha_0(\vect{x}) \Should y^{\text{nom}} \label{eq:y_0} \\\label{eq:y_1}  \tag{4.2}
        \dot{y} &=  \frac{\partial \alpha_0(\vect{x})}{\partial\vect{x}} \dot{\vect{x}} = \frac{\partial \alpha_0(\vect{x})}{\partial\vect{x}}\vect{f}_1 (\vect{x}) + \underbrace{\frac{\partial \alpha_0(\vect{x})}{\partial\vect{x}}\vect{f}_{2,1} (\vect{x})}_{=0}u + \frac{\partial \alpha_0(\vect{x})}{\partial\vect{x}}\vect{f}_{2,2} (\vect{x})\rho
        \\ \nonumber & \coloneqq   \alpha_1(\vect{x},\rho) \Should 0 \\  \tag{4.3}
        y^{(2)} &=  \frac{\partial \alpha_1(\vect{x},\rho)}{\partial\vect{x}} \dot{\vect{x}} + \frac{\partial \alpha_1(\vect{x},\rho)}{\partial\rho} \dot{\rho} \coloneqq   \alpha_2(\vect{x},\boldsymbol{\varphi}) \Should 0 \\
        &~~~~~~~~~~~~~~\vdots \nonumber 
    \end{align}    
    \begin{align}
        \label{eq:y_n}  \tag{4.n+1}
        y^{(n)} &=  
        \frac{\partial \alpha_{n-1}(\vect{x},\boldsymbol{\varphi})}{\partial\vect{x}}\underbrace{\left( \vect{f}_1 (\vect{x}) + \vect{f}_{2,1} (\vect{x})u + \vect{f}_{2,2} (\vect{x})\rho \right)}_{\dot{\vect{x}}} 
        \\ \nonumber & + \frac{\partial\alpha_{n-1}(\vect{x},\boldsymbol{\varphi})}{\partial \underbrace{\left(  \rho,   \dots  ,\rho^{(\delta-2)},  \rho^{(\delta-1)}\right)^T}_{\boldsymbol{\varphi}}} \underbrace{\left( \begin{array}{c}  \dot{\rho} \\\vdots\\ \rho^{(\delta-1)} \\ \rho^{(\delta)} \end{array} \right)}_{\dot{\boldsymbol{\varphi}}} 
        \\ \nonumber
        & = \underbrace{\frac{\partial\alpha_{n-1}(\vect{x},\boldsymbol{\varphi})}{\partial\vect{x}}\left( \vect{f}_1 (\vect{x}) + \vect{f}_{2,2}(\vect{x})\rho \right) + \frac{\partial\alpha_{n-1}(\vect{x},\boldsymbol{\varphi})}{\partial \left(  \rho,   \dots  ,\rho^{(\delta-2)} \right)^T} \left( \begin{array}{c}  \dot{\rho} \\\vdots\\ \rho^{(\delta-1)} \end{array} \right)}_{\coloneqq \alpha_n(\vect{x},\boldsymbol{\varphi})} 
        \\ \nonumber & + \underbrace{\frac{\partial\alpha_{n-1}(\vect{x},\boldsymbol{\varphi})}{\partial\vect{x}} \vect{f}_{2,1} (\vect{x})}_{\coloneqq \beta_u (\vect{x},\boldsymbol{\varphi})} u + \underbrace{\frac{\partial \alpha_{n-1}(\vect{x},\boldsymbol{\varphi})}{\partial\rho^{(\delta-1)}}}_{\coloneqq \beta_{\rho}(\vect{x},\boldsymbol{\varphi})} \rho^{(\delta)}
        \\ \nonumber 
        &=  \alpha_n(\vect{x},\boldsymbol{\varphi}) + \beta_u (\vect{x},\boldsymbol{\varphi}) u + \beta_{\rho}(\vect{x},\boldsymbol{\varphi}) \rho^{(\delta)} \Should 0 ~~
        \\ \nonumber & \text{with nonlinear functions } \alpha_n(\vect{x}, \boldsymbol{\varphi}) \text{, }\beta_u (\vect{x},\rho) \text{, and } \beta_{\rho} (\vect{x},\boldsymbol{\varphi})
    \end{align}
\addtocounter{equation}{1}
Note that the number of necessary differentiations will typically be small as processes are usually designed following the principle of local disturbance rejection \citep{Skogestad.2001}.
According to local disturbance rejection, control inputs $u$ should affect controlled outputs $y$ as directly as possible, i.e., small relative degrees $r$ are generally preferred.
Equation (\ref{eq:y_n}) shows that the derivative $y^{(n)}$ is influenced by the scheduling decision, i.e., the $\delta$-th derivative of the production rate $\rho^{(\delta)}$, and by the control decision, i.e., the input $u$.
The integer $\delta$ defines the order for our dynamic ramping constraints and thus the ramping degree of freedom $\nu$, which acts as free variable in scheduling optimization, must be equal to $\rho^{(\delta)}$.
In other words: Control can only hold the output $y$ constant if scheduling ramps the production rate with a $\delta$-th order dynamic.
If a ramping dynamic with order $\gamma < \delta$ would be chosen, the $\gamma$-th derivative would be the ramping degree of freedom, i.e., scheduling optimization could perform a step-change on $\rho^{(\gamma)}$.
However, this step-change would act on a derivative $y^{(k)}$ with $k<n$, a derivative that is not influenced by the control input $u$.
Consequently, control has no handle to hold $y^{(k)}$ at zero and thus has no handle to hold $y$ constant.  
Thus, ramping with an order smaller than $\delta$ must lead to a deviation of the output $y$ from its nominal value irrespective of the used control method.

For the ramping limits $\nu^{\text{min}}$ and $\nu^{\text{max}}$, we rearrange Equation~(\ref{eq:y_n}) to get the ramping degree of freedom $\nu$, which is equal to the highest derivative $\rho^{(\delta)}$, as function of the input $u$:
    \begin{align} \label{eq:nu}
        \nu = \frac{-\alpha_n(\vect{x},\boldsymbol{\varphi}) - \beta_u (\vect{x},\boldsymbol{\varphi}) u}{\beta_{\rho}(\vect{x},\boldsymbol{\varphi})}
    \end{align}
In Equation (\ref{eq:nu}), the ramping degree of freedom $\nu$ can be influenced by the input $u$, i.e., the bounds on $u$ limit $\nu$. 
If $\beta_u (\vect{x},\boldsymbol{\varphi})$ and $\beta_{\rho} (\vect{x},\boldsymbol{\varphi})$ have the same sign, we get the bounds by:
    \begin{align}
        \label{eq:nu_min}
        \nu^{\text{min}}(\vect{x},\boldsymbol{\varphi})  = \frac{-\alpha_n(\vect{x},\boldsymbol{\varphi}) - \beta_u (\vect{x},\boldsymbol{\varphi}) u^{\text{max}}}{\beta_{\rho}(\vect{x},\boldsymbol{\varphi})} \\
        \label{eq:nu_max}
        \nu^{\text{max}}(\vect{x},\boldsymbol{\varphi})  = \frac{-\alpha_n(\vect{x},\boldsymbol{\varphi}) - \beta_u (\vect{x},\boldsymbol{\varphi}) u^{\text{min}}}{\beta_{\rho}(\vect{x},\boldsymbol{\varphi})}
    \end{align}
If $\beta_u (\vect{x},\boldsymbol{\varphi})$ and $\beta_{\rho} (\vect{x},\boldsymbol{\varphi})$ have a different sign, $u^{\text{max}}$ and $u^{\text{min}}$ need to be swapped in Equations (\ref{eq:nu_min}) and (\ref{eq:nu_max}). 
From Equation (\ref{eq:nu_min}) follows that if scheduling chooses the ramping degree of freedom to be smaller than $\nu^{\text{min}}(\vect{x},\boldsymbol{\varphi})$, control cannot hold the output $y$ constant because for $\nu < \nu^{\text{min}}(\vect{x},\boldsymbol{\varphi})$ a control input $u > u^{\text{max}}$ would be needed. 

As the limits $\nu^{\text{min}}(\vect{x},\boldsymbol{\varphi})$ and $\nu^{\text{max}}(\vect{x},\boldsymbol{\varphi})$ depend on the state vector $\vect{x}$, we want to express $\vect{x}$ as a function of the ramping state vector $\boldsymbol{\varphi}$ such that we can express the limits purely based on $\boldsymbol{\varphi}$.
To this end, we make use of the assumption that the control maintains $y$ at its nominal value $y^{\text{nom}}$ and receive the following system of equations (compare to Equations (\ref{eq:y_0}) - (4.n):
    \begin{align}
        \label{eq:system_gamma}
        \vect{0} = 
        \underbrace{\left(\begin{array}{c} y^{\text{nom}} \\ 0 \\\vdots \\0\end{array} \right) -  \left(\begin{array}{c} \alpha_0(\vect{x}) \\ \alpha_1(\vect{x},\boldsymbol{\varphi}) \\ \vdots \\ \alpha_{n-1}(\vect{x},\boldsymbol{\varphi}) \end{array} \right)}_{\boldsymbol{\alpha}(\vect{x},\boldsymbol{\varphi})} 
    \end{align}
Equation~\ref{eq:system_gamma} gives a system of $n$ nonlinear equations which implicitly define the state vector $\vect{x}$ as a function of the ramping state vector $\boldsymbol{\varphi}$.
According to the implicit function theorem, $\vect{x}$ can be calculated as a
locally unique function $\boldsymbol{\Gamma}(\boldsymbol{\varphi})$ if the $n\times n$ Jacobian matrix $\vect{J}(\vect{x},\boldsymbol{\varphi}) =\left(\frac{\partial\boldsymbol{\alpha}(\vect{x},\boldsymbol{\varphi})}{\partial\vect{x}} \right)$ has a non-zero determinant. 
To check if a unique function $\boldsymbol{\Gamma}(\boldsymbol{\varphi})$ exists over the complete operating range, we analytically evaluate the determinant using computer algebra and check if the determinant is nonzero over the complete operating range.
If symbolically checking whether the determinant of the the Jacobian is always nonzero is not possible, it might be reasonable to evaluate the determinant of the Jacobian for a specific point, e.g., the nominal operating point, and then proceed with trying to solve the equation system.
Note that the existence of a function $\boldsymbol{\Gamma}$ that gives the original states $\vect{x}$ is a condition for exact linearization and controllability.
In other words, for a well-designed process, where the output of interest $y$ is, in fact, controllable using the input $u$ over the complete operating range, the determinant of the Jacobian is nonzero.

In this work, we calculate $\boldsymbol{\Gamma}(\boldsymbol{\varphi})$ symbolically using the computer algebra package SymPy \citep{meurer2017sympy}.
With $\boldsymbol{\Gamma}(\boldsymbol{\varphi})$, we can calculate the limits $\nu^{\text{min}}$ and $\nu^{\text{max}}$ as functions of the ramping state vector $\boldsymbol{\varphi}$ and receive the high-order dynamic ramping constraints (Equations \ref{eq:dyn_ramping_a} and \ref{eq:dyn_ramping_b}).

Note that the assumptions necessary to derive the dynamic ramping constraints are very restrictive and limit their applicability.
In Section \ref{sec:discussion}, we discuss possible extensions.

\subsection{Piecewise affine limits}
\label{subsec:PWL_limits}
\begin{figure}[h!]
\centering
\footnotesize
\includegraphics{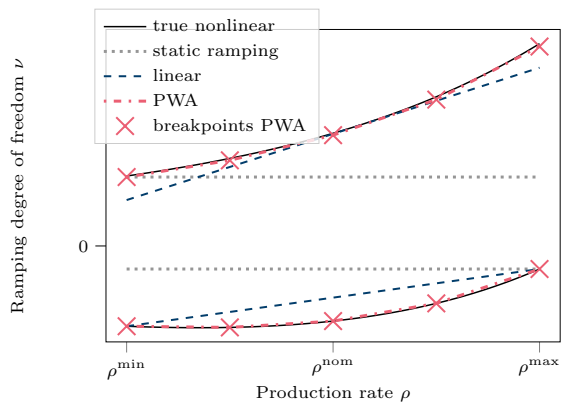}
    \caption{Constraints for ramping degree of freedom $\nu$ as function of production rate $\rho$. True nonlinear limits in comparison to static ramping limits, linear limits, and piece-wise affine (PWA) limits for an illustrative case with first-order dynamics. For first-order dynamics, the transformed state vector $\boldsymbol{\varphi}$ is of dimension one and equal to the production rate $\rho$. Consequently, the limits on the ramping degree of freedom $\nu$ only depend on $\rho$.}
    \label{fig:comp_bounds}
\end{figure}

In general, the limits $\nu^{\text{min}}(\boldsymbol{\varphi})$, $\nu^{\text{max}}(\boldsymbol{\varphi})$ (Equations \ref{eq:nu_min} and \ref{eq:nu_max}) are nonlinear functions and linear or piece-wise affine approximations are needed to achieve an MILP formulation.
Fortunately, if both limits are approximated conservatively, the resulting trajectory is always feasible because the ramping degree of freedom $\nu$ stays within the true feasible range.
In contrast to static ramping constraints, approximating the nonlinear limits allows to balance conservativeness against computational complexity as we discuss in the following for the case of first-order dynamics.
If static ramping constraints or purely linear functions are chosen, computational costs are small.
However, potentially, a high conservatism has to be accepted as a large amount of the feasible region for the ramping degree of freedom $\nu$ is cut off. 
In Figure \ref{fig:comp_bounds}, the feasible region between the nonlinear limits is much bigger than the feasible region between the static limits and also bigger than the feasible region between the linear limits. 
Using piece-wise affine functions, the conservative limits can enclose more of the feasible nonlinear region. 
However, if the feasible region is non-convex, binary variables have to be introduced, which significantly increase the computational burden in optimization.

In the case of first-order dynamics, the true limits on the ramping degree of freedom $\nu$ can simply be plotted to choose appropriate bounds.
In case of high-order dynamics, the limits are multivariate functions and multivariate regression methods that give piece-wise affine functions, e.g., hinging hyperplanes \citep{L.Breiman.1993,A.A.Adeniran.2017,Kamper.2021b}, convex region surrogates \citep{Zhang.2016b,Schweidtmann.2021}, or artificial neural networks with ReLU activation functions \citep{BjarneGrimstad.2019}, must be used. 
However, as discussed after Equation~(4n+1), typically the number of necessary differentiations is expected to be small. 

By employing a conservative approximation of the true nonlinear ramping limits, the feasibility of the found solution on the original nonlinear model is guaranteed.
As feasible area is cut off, the found optimum of the approximated MILP problem might deviate from the optimal solution of the original MINLP problem, i.e., optimality on the original problem might not be achieved.
However, with piecewise affine functions, the true nonlinear limits can, in principle, be approximated to any accuracy if a sufficient number of piecewise linear elements is used.
These piecewise linear elements however can increase the number of binary variables leading to higher computational cost.
Thus, a suitable trade-off between accuracy of the approximation and number of binary variables must be found.

\section{Case Study}
\label{sec:case_study}
As case study, we consider two continuous stirred tank reactors (CSTRs) with different dynamic orders in parallel configuration.
The first CSTR is a dimensionless benchmark CSTR from literature \citep{FloresTlacuahuac.2006,Du.2015} described by material and energy balances. 
Note that, for this CSTR, the material flow rate equals the production rate.
As we use the symbol $\rho$ for the production rate throughout the paper, we denote the material flow rate as $\rho$ to preserve consistency, even if this is an unusual choice.
    \begin{align}
        \label{eq:reactor_1a} \tag{CSTR1a}
        &\dot{c} =  (1 - c)\frac{\rho}{V}- c k e^{- \frac{N}{T}} \\
        \label{eq:reactor_1b} \tag{CSTR1b}
        &\dot{T} = (T_{f} - T)\frac{\rho}{V} + c k e^{- \frac{N}{T}} - F_{c} \alpha_c \left(T - T_{c}\right)
    \end{align}
The CSTR states $\vect{x}$ are the concentration $c$ and the temperature $T$. 
The control input $u$ is the coolant flow rate $F_c$ with bounds $F_c^{min} = 0 \frac{1}{h}$, $F_c^{max} = 700 \frac{1}{h}$ \citep{Du.2015}. The material flow rate $\rho$ is a degree of freedom and equals the production rate. 
All other symbols are constant parameters given in Table~\ref{tab:cstr_parameter}.
Following \cite{Schafer.2020}, we assume an oversizing of 20~\% such that the dimensionless production rate can be varied by +/- 20~\% around the nominal value $\rho^{nom} = 1\frac{1}{h}$, i.e., $\rho^{\text{min}} = 0.8\frac{1}{h}$ and $\rho^{\text{max}} = 1.2\frac{1}{h}$.
Additionally, we assume that the concentration has to be maintained constant at $c=c^{nom}=0.1367$ (product 2 in \cite{FloresTlacuahuac.2006}). Thus, the output $y$ is given  by $y=h(\vect{x})=c$.
The second CSTR is identical to the first CSTR except for an additional inertia in the form of a cooling jacket with temperature $T_j$.
With an energy balance of the cooling jacket, the model of CSTR 2 reads:
    \begin{align}
        \label{eq:reactor_2a} \tag{CSTR2a}
        & \text{Equation (\ref{eq:reactor_1a})}\\
        \label{eq:reactor_2b} \tag{CSTR2b}
        &\dot{T} = (T_{f} - T)\frac{\rho}{V} + c k e^{- \frac{N}{T}} + \tau_1(T_j - T) \\ \label{eq:reactor_2c}\tag{CSTR2c}
        &\dot{T}_j = \tau_2(T - T_j) - F_{c} \alpha_c \left(T_j - T_{c}\right)
    \end{align}
The inverse time constants $\tau_1$, $\tau_2$ are derived from the reactor studied in \cite{M.Mezghani.2002} and given in Table~\ref{tab:cstr_parameter}.
\begin{table}[bt]
    \centering
    \caption{Dimensionless CSTR model parameters from \cite{FloresTlacuahuac.2006,Du.2015} and time constants from \cite{M.Mezghani.2002}.}
    \label{tab:cstr_parameter}
    \begin{tabular}{llr}
     &symbol & value  \\
     \hline
    volume & $V$ & $20$  \\
    reaction constant & $k$ & $300\frac{1}{h}$ \\
    activation energy & $N$ & 5 \\
    feed temperature & $T_f$ & 0.3947 \\
    heat transfer coefficient & $\alpha_c$ & $1.95\cdot 10^{-4}$ \\
    coolant temperature & $T_c$ & 0.3816 \\
    inverse time constant & $\tau_1$ & $4.84\frac{1}{h}$ \\
    inverse time constant & $\tau_2$ & $14.66\frac{1}{h}$ \\
    \hline  
    \end{tabular}
\end{table}
    
\begin{figure}[h]
    \centering
 \includegraphics{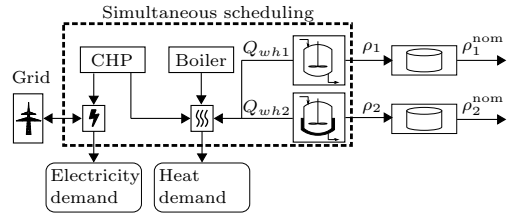}
  \caption{Case Study: Simultaneous scheduling of two continuous stirred tank reactors (CSTRs) with variable production rates $\rho_1$ and $\rho_2$, a boiler, and a combined heat and power plant (CHP). The waste heat of the two CSTRs, $Q_{wh1}$ and $Q_{wh2}$, is used to partly satisfy a non-flexible heat demand. Additionally, a non-flexible electricity demand has to be fulfilled.}
   \label{fig:case_study}
\end{figure}

We study the two CSTRs in combination with a benchmark multi-energy system from \cite{SusanneSass.2020}. The waste heat of the two flexible CSTRs is integrated to partly satisfy the heat demand of inflexible consumers (Figure \ref{fig:case_study}). In other words: The two CSTRs demand cooling but the removed heat is integrated such that the CSTRs supply waste heat. The dimensionless waste heat $Q_{wh}$ is:
\begin{align}
    Q_{wh} = F_{c} \alpha_c \left(T - T_{c}\right)
\end{align}

The multi-energy system consists of a combined heat and power plant (CHP) and a boiler (Bo) that together need to satisfy a heat demand and an electricity demand. 
Electricity can be bought from the grid and sold to the grid.
As the demands in \cite{SusanneSass.2020} (available online at \cite{HECI}) are calculated using weather data from the 28th November  2018, we use  the German day-ahead market price series for electricity from the same day \citep{SMARD}.
We scale the waste heat of the dimensionless benchmark CSTRs such that the nominal waste heat of each CSTR corresponds to 10\% of the maximum heat demand. 
Both CHP and boiler have  a variable part-load efficiency and a minimum part-load \citep{SusanneSass.2020}.

\subsection{Derivation of dynamic ramping constraints}
\begin{figure}[ht]   
\centering
\footnotesize
\includegraphics{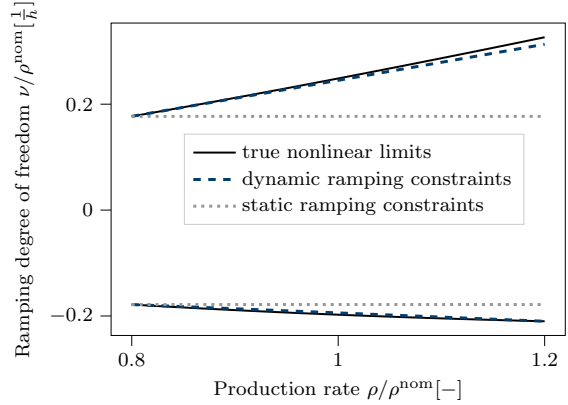}
    \caption{True nonlinear limits and linear dynamic ramping constraints for the non-jacketed CSTR 1. Static ramping constraints are indicated for comparison.}
    \label{fig:bounds}
\end{figure}

The first CSTR has two differential states and is exactly input-state linearizable, i.e., the relative degree $r$ is 2. 
Thus, the output $y$ has to be differentiated two times until the input $u=F_c$ appears explicitly:
    \begin{align}
        \label{eq:CSTR1_h}
        y=&c \coloneqq \alpha_0(c) \Should c^{\text{nom}} 
        \\
        \label{eq:CSTR1_diff1}
        \dot{y} =&  \frac{\partial \alpha_0(c)}{\partial c} \dot{c} = (1 - c)\frac{\rho}{V}- c k e^{- \frac{N}{T}} \coloneqq \alpha_1(c,T,\rho)  \Should 0  \\
        \label{eq:CSTR1_diff2}
        y^{(2)} =& \frac{\partial \alpha_1(c,T,\rho)}{\partial c} \dot{c} + \frac{\partial \alpha_1(c,T,\rho)}{\partial T} \dot{T}+ \frac{\partial \alpha_1(c,T,\rho)}{\partial \rho} \dot{\rho}  \\ \nonumber
        =& - \left[\frac{\rho}{V} + k e^{- \frac{N}{T}}\right]\left[(1 - c)\frac{\rho}{V}
        - c k e^{- \frac{N}{T}}\right] 
        \\ \nonumber & - \left[  \frac{c k N e^{- \frac{N}{T}}}{T^2}\right]\left[(T_{f} - T)\frac{\rho}{V} + c k e^{- \frac{N}{T}} - F_{c} \alpha_c \left(T - T_{c}\right) \right]  + \left[\frac{1-c}{V} \right]\dot{\rho} \\ \nonumber
        =& \alpha_2(\vect{x},\rho) + \beta_u (\vect{x},\rho) F_c + \beta_{\rho}(\vect{x},\rho)\nu \Should 0 ,
        \\ \nonumber
        &\text{with } \nu = \dot{\rho},~ \alpha_2(\vect{x},\rho) = - \left[\frac{\rho}{V} + k e^{- \frac{N}{T}}\right]\left[(1 - c)\frac{\rho}{V}- c k e^{- \frac{N}{T}}\right]  
        \\ \nonumber & - \left[  \frac{c k N e^{- \frac{N}{T}}}{T^2}\right]\left[(T_{f} - T)\frac{\rho}{V} + c k e^{- \frac{N}{T}} \right] ,
        \\& \nonumber\beta_u (\vect{x},\rho) =   \frac{c k N e^{- \frac{N}{T}} \alpha_c \left(T - T_{c}\right)}{T^2} > 0, \text{ and }  \beta_{\rho}(\vect{x},\rho) = \frac{1-c}{V}  > 0 
    \end{align}
Consequently, the production rate $\rho$ can be changed with a first-order dynamic because only the first derivative of the production rate appears during the differentiations.
The equation system to calculate the transformation $\vect{x}=\boldsymbol{\Gamma}(\rho)$ (compare to Equation (\ref{eq:system_gamma})) is given by Equations (\ref{eq:CSTR1_h}) and (\ref{eq:CSTR1_diff1}).
The determinant of the Jacobian is $-  \frac{c k N e^{- \frac{N}{T}}}{T^2}$ (compare to discussion after Equation~(\ref{eq:system_gamma})), which is always nonzero as the states $c, T$ are always nonzero in the considered operating range, and the parameters $k, N$ are also nonzero.
Consequently, the two equations can be solved to calculate the state vector $\vect{x}$ as function of the production rate $\rho$:
    \begin{align}
        \underbrace{\left(\begin{array}{c} c \\ T\end{array} \right)}_{\vect{x}} =\underbrace{\left(\begin{array}{c} c^{nom} \\ \frac{N}{\ln{\left( \frac{V c k}{\rho \left(1 - c\right)} \right)}}\end{array} \right)}_{\boldsymbol{\Gamma}(\rho)} 
    \end{align}
With $\boldsymbol{\Gamma}(\rho)$, the limits of the ramping degree of freedom $\nu$ are calculated as a function of the production rate $\rho$ using Equations (\ref{eq:nu_min}) and (\ref{eq:nu_max}) (Figure \ref{fig:bounds}).
Based on the visualization in Figure \ref{fig:bounds}, we choose purely linear dynamic ramping limits $\nu^{\text{min}}(\rho)$, $\nu^{\text{max}}(\rho)$ as these give a good approximation.
Note that in many control engineering applications, linear approximations give a reasonable approximation of nonlinear functions as long as the process is close to a nominal operating point \citep{Corriou.2018} like in the example here.

\begin{figure}[ht] 
\centering
\footnotesize
\includegraphics{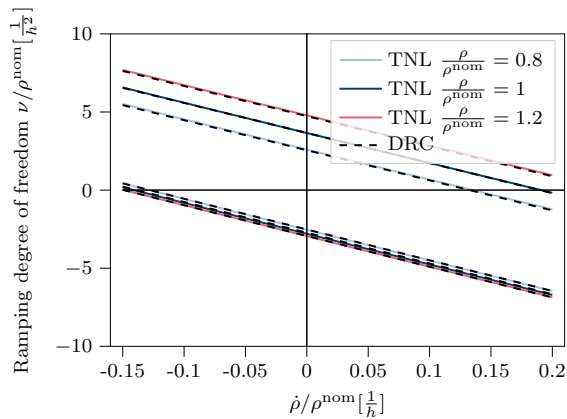}
    \caption{True nonlinear limits (TNL) and corresponding linear dynamic ramping constraints (DRC) for the jacketed CSTR 2. The limits are functions of the production rate $\rho$ and its first derivative $\dot{\rho}$. In this figure, the limits are given depending on $\dot{\rho}$ for 3 different values of $\rho$.}
    \label{fig:bounds_rwj}
\end{figure}

For the jacketed CSTR 2, the output $c$ must be differentiated three times until the input $u$ appears.
For this CSTR, we only show the functional dependencies 
to avoid lengthy equations:
    \begin{align}
        \label{eq:CSTR2_h}
        y&=c \coloneqq \alpha_0(c) \Should c^{\text{nom}} \\
        \label{eq:CSTR2_diff1}
        \dot{y} &=  \frac{\partial \alpha_0(c)}{\partial c} \dot{c}\coloneqq \alpha_1(c,T,\rho) \Should 0   \\
        y^{(2)} &= \frac{\partial \alpha_1(c,T,\rho)}{\partial c} \dot{c} + \frac{\partial \alpha_1(c,T,\rho)}{\partial T} \dot{T} + \frac{\partial \alpha_1(c,T,\rho)}{\partial \rho} \dot{\rho} \coloneqq \alpha_2(c,T,T_j,\rho,\dot{\rho}) \Should 0 \\
        \label{eq:CSTR2_y3}
        y^{(3)} &= \frac{\partial \alpha_2(c,T,T_j,\rho,\dot{\rho})}{\partial \vect{x}} \left(\begin{array}{c} \dot{c} \\ \dot{T} \\ \dot{T}_j \end{array} \right) + \frac{\partial \alpha_2(c,T,T_j,\rho,\dot{\rho})}{\partial \boldsymbol{\varphi}} \left(\begin{array}{c} \dot{\rho} \\ \rho^{(2)}  \end{array} \right) \\ \nonumber
        &= \alpha_3(\vect{x},\boldsymbol{\varphi}) + \beta_u (\vect{x},\boldsymbol{\varphi}) F_c + \beta_{\rho}(\vect{x},\boldsymbol{\varphi})\nu \Should 0 \text{ with } \nu = \rho^{(2)}
    \end{align}
The complete derivatives are given in the Supplementary Information.
For CSTR 2, the Jacobian of Equation~(\ref{eq:system_gamma}) is $\frac{ N^{2} \alpha_{c} c^{2} k^{2} e^{- \frac{2  N}{T}}}{T^{4}}$, which is always nonzero. 
Thus, the states $\vect{x}$ can be given as a function of the ramping state vector  $\boldsymbol{\Gamma}(\boldsymbol{\varphi})$.
As the second time derivative of the production rate, $\rho^{(2)}$, appears in $y^{(3)}$, this second derivative $\rho^{(2)}$ is the ramping degree of freedom $\nu$ and the ramping state vector $\boldsymbol{\varphi}$ is two-dimensional
    $\boldsymbol{\varphi} = \left(\begin{array}{c} \rho \\ \dot{\rho} \end{array} \right)$.
In other words: A second-order ramping constraint with $\delta=2$ in Equation \ref{eq:dyn_ramping_a} is needed.
We choose the simplified ramping limits $\nu^{\text{min}}(\rho, \dot{\rho} )$, $\nu^{\text{max}}(\rho, \dot{\rho} )$ to be linear in both $\rho$ and $\dot{\rho}$:
    \begin{align}
        &\nu^{\text{min}}(\rho, \dot{\rho} ) = \nu_0^{\text{min}} + m_{\rho}^{\text{min}}\rho + m_{\dot{\rho}}^{\text{min}}\dot{\rho} \\
        &\nu^{\text{max}}(\rho, \dot{\rho} ) = \nu_0^{\text{max}} + m_{\rho}^{\text{max}}\rho + m_{\dot{\rho}}^{\text{max}}\dot{\rho}, 
    \end{align}
    with parameters $\nu_0^{\text{min}}$, $m_{\rho}^{\text{min}}$, $m_{\dot{\rho}}^{\text{min}}$, $\nu_0^{\text{max}}$, $m_{\rho}^{\text{max}}$, $m_{\dot{\rho}}^{\text{max}}$.
To parameterize the bounds, we sample the operating range using 100 equally distributed points for $\rho$ and 100 equally distributed points for $\dot{\rho}$ such that there are 10,000 points in total.
For each point, the true nonlinear limits are calculated, and the parameters are fitted to the nonlinear limits using the normal equation method \citep{Lewis.2006}. 
To make the resulting limits conservative, we search through the grid for the highest violation of the nonlinear limits and adapt $\nu_0^{\text{min}}$ and $\nu_0^{\text{max}}$ to the safe side.
As visualized in Figure \ref{fig:bounds_rwj}, again the linear limits are close to the nonlinear limits.
This finding could be expected because linear limits already give a reasonable approximation for the non-jacketed CSTR1 and the additional heat transfer terms in Equations (\ref{eq:reactor_2b}) and (\ref{eq:reactor_2c}) are purely linear.

\begin{figure*}[ht] 
\centering
\footnotesize
\includegraphics{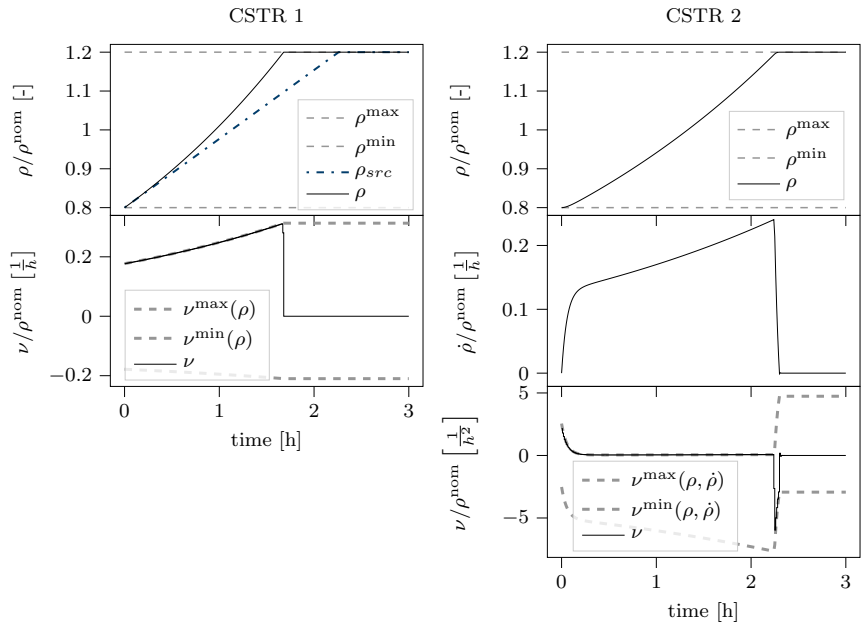}
    \caption{Fastest possible ramp up from minimum production rate $\rho^{\text{min}}$ to maximum production rate $\rho^{\text{max}}$ for CSTR 1 and CSTR 2. The ramping degree of freedom $\nu$ is shown together with its dynamic limits $\nu^{\text{min}}(\cdot)$, $\nu^{\text{max}}(\cdot)$. For CSTR 1, the fastest possible ramp up with a static ramping constraint $\rho_{src}$ is shown for comparison.}
    \label{fig:ramp_up}
\end{figure*}

To illustrate the dynamic ramping constraints derived for the non-jacketed CSTR 1 and the jacketed CSTR 2, we perform a first optimization in which we ramp the reactors from minimum production rate to maximum production rate as fast as possible (Figure \ref{fig:ramp_up}). 
For CSTR 1, the ramp up takes 1.7~h while it would take 2.3~h (+35~\%) with a static ramping constraint where the maximum ramping $\nu^{\text{max}}$ is constant (Figure \ref{fig:bounds}).
The ramp up of the jacketed CSTR 2 takes more than half an hour longer compared to CSTR 1 because of the cooling jacket inertia.
We cannot compare the ramp-up of CSTR 2 to traditional first-order static ramping because $y^{(3)}$ can only be held constant by manipulating $F_c$ if the second derivative of the production rate $\rho^{(2)}$ is defined (compare to Equation~(\ref{eq:CSTR2_y3})). 
A first-order ramp would give a step change on the first derivative $\dot{\rho}$ such that $\rho^{(2)}$ is not defined and the output $y$ must deviate from the nominal value. 
Depending on the application, such deviations might be acceptable, still, they can be avoided using a second-order ramp. 
Moreover, if deviations are only acceptable up to a certain tolerance, it might be necessary to choose a first order ramp slow enough such that deviations can be corrected by the underlying control while a second-order ramp might allow a faster ramp-up.

As linear dynamic ramping constraints are used for the two CSTRs, we expect the resulting formulation to be computationally efficient. 
Moreover, the two shortcomings of static ramping constraints can be illustrated: While the jacketed CSTR2 cannot be modeled accurately by static ramping constraints due to second-order dynamics, dynamic ramping constraints ramp up the non-jacketed CSTR1 35~\% faster than static ramping constraints.
\subsection{Waste heat model}
For the dynamic MILP scheduling problem P2, a process energy demand model is needed (\ref{eq:P2_process_energy_dem}).
In this case study, the energy demand of the CSTRs corresponds to the waste heat removed from the CSTRs (cf. Figure \ref{fig:case_study}).
For the scheduling waste heat model, we start by investigating purely linear functions as approximation
    \begin{align}
        &Q_{wh1} = a_{0,1} + a_{1,1}\rho_1 + a_{2,1}\nu_1 , \\
        &Q_{wh2} = a_{0,2} + a_{1,2}\rho_2 + a_{2,2}\dot{\rho}_2 + a_{3,2}\nu_2 ,
    \end{align}
    with coefficients $a_{0,1}$, $a_{1,1}$, $a_{2,1}$, $a_{0,2}$, $a_{1,2}$, $a_{2,2}$, $a_{3,2}$.
To determine the coefficients, we sample the operating region using 11 equally distributed points for $\rho$, $\nu$, and $\dot{\rho}$ (only CSTR2), which gives $11\times11$ points for CSTR 1  and $11\times11\times11$ points for CSTR 2.
For each point, we calculate the waste heat from the nonlinear model and subsequently fit the coefficients using the normal equation method \citep{Lewis.2006}.
The average absolute deviation between fit and nonlinear model is 4 \% of the nominal waste heat for CSTR 1 and 1 \% of the nominal waste heat for CSTR 2. 
As these deviations are small, we do not study potentially more accurate piece-wise affine models for the waste heat.
Moreover, we assume that small deviations can be compensated by the energy system components, as these components typically react much faster than the chemical process.

\subsection{Optimization problem}
Based on the dynamic ramping constraints and the waste heat model, in this subsection, we formulate the scheduling optimization problem P2.

Following \cite{Voll.2014}, the efficiency curves (Equation (\ref{eq:P1_eff})) for CHP and boiler specified in \cite{SusanneSass.2020} are discretized with one affine element to obtain a satisfactory discretization.
Note that a small number of piece-wise affine elements, i.e., one or two, is often sufficient for modeling univariate efficiency functions of typical energy system components \citep{AndreasKamper.2021}.

For the storage units, we enforce that the filling level is 50~\% of the maximum filling level at the beginning and the end of the day such that the total production is fixed. 
We assume that the capacity of the two storage units is equal to 3 hours of production at nominal production rate.
A storage capacity of 3 hours  is the lower bound of the range studied by \cite{Schafer.2020} who consider storage capacities between 3 and 48 hours for demand response applications.
We find that the storage capacity is no limiting factor in our case study, i.e., the bounds of the storage filling level are never reached. Additionally, the optimization is repeated without storage constraints and the result does not change.

As electricity prices and energy demands are given with one hour time resolution, we use one hour timesteps for the ramping degrees of freedom $\nu_1$, $\nu_2$, and for the on/off binaries $z_i^{\text{on}}$ (cf. Equation \ref{eq:P1_part_load}).
The remaining variables are discretized with 2 collocation elements per hour and 4 points per element.
Note that we calculate economic profit in a simulation and thereby automatically verify the adequacy of our time discretization.

The optimization problem is formulated using pyomo \citep{hart2017pyomo,hart2011pyomo} and pyomo.dae \citep{Nicholson.2018} for discretization.
We solve the optimization problem using gurobi version 8.1.0 \citep{gurobi} with zero optimality gap.
All calculations are performed on a Windows 10 machine with an Intel(R) Core(TM) i5-8250U core and 24 GB RAM.
\subsection{Results}
The cost reductions achieved by waste heat integration relate to an operation without waste heat integration.
First, we calculate the cost reductions through waste heat integration in a scenario without DR, i.e., both CSTRs are operated in steady state such that waste heat production is constant.
Second, in the DR scenario, we optimize the operation of energy system and CSTRs simultaneously, i.e., we solve the dynamic ramping problem (P2). 
Next, we simulate the resulting schedule on the original nonlinear process model over the 24 hour time horizon to determine the cost reduction achieved through waste heat integration.
Third, to bound the DR potential, we solve the original MINLP problem (P1) using the solver BARON version 20.10.16 \citep{Khajavirad.2018}.
As BARON with default settings does not provide a feasible point after 2 hours of optimization, we generate a point by fixing the integer on/off-decisions to the values from the solution of our dynamic ramping problem (P2), i.e., we reduce the MINLP to a nonlinear program (NLP).
For this NLP, BARON provides a feasible point, which we use as initial point to solve the original MINLP (P1) with BARON.
After 2 hours of calculation, we use the reported lower bound to determine the maximum possible DR potential. 

The DR scenario reduces 40.8~\% more costs than the steady-state scenario (Table \ref{tab:results}).
This DR improvement is only slightly below the MINLP bound, which is 42.7~\%.
Further, the optimization runtime of the dynamic ramping problem (P2) is below 2 seconds, which is short enough for online application.

\begin{table*}[h]
    \caption{Cost reduction for the considered day achieved by waste heat integration for different scenarios. For the MINLP solution the lower bound found after 2 hours of optimization is reported.}
    \label{tab:results}
    \begin{tabular}{lrr}
     Scenario & Cost reduction [\euro] & DR improvement [\%] \\
    \hline
    No DR (steady-state, operation) &  90.11 & -  \\
    DR with dynamic ramping (P2) &  126.91 & 40.8\\
    MINLP bound (P1) & 128.61  & 42.7\\
    \hline  
    \end{tabular}
\end{table*}

\begin{figure}[h!]        
\centering
\includegraphics{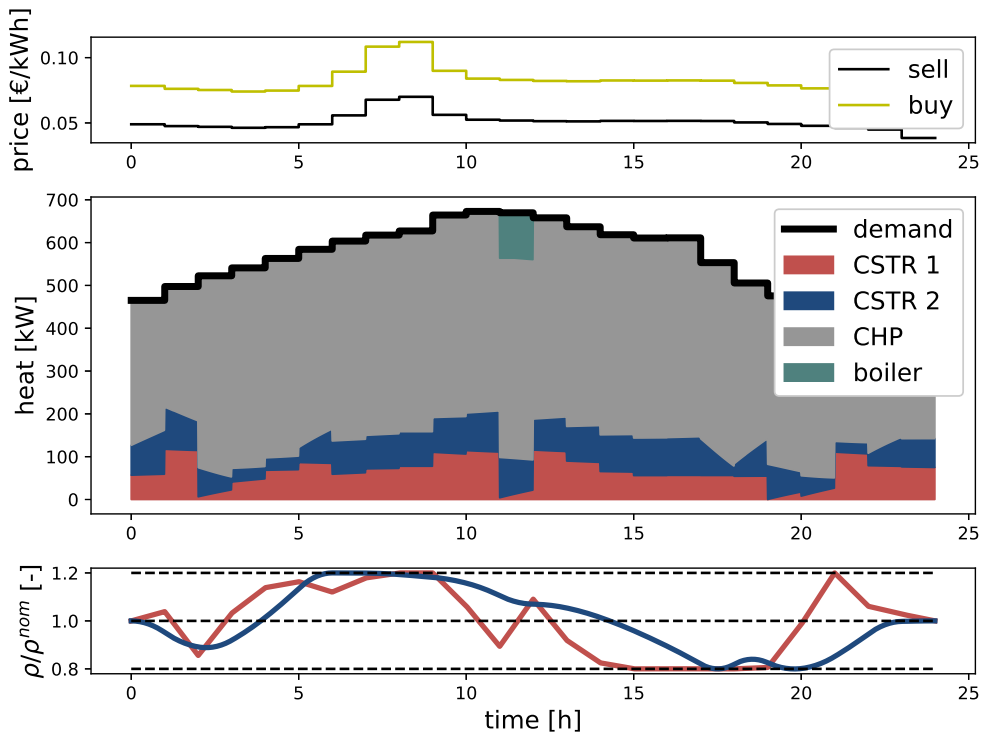}
    \caption{Resulting operation in the DR scenario with dynamic ramping constraints showing the heat supplied by the two CSTRs, the CHP, and the boiler, and the production $\rho$ delivered by the two CSTRs (in red for CSTR 1 and blue for CSTR 2). }
    \label{fig:result}
\end{figure}

Interestingly, in the DR scenario, the boiler is only active for 1 hour (Figure \ref{fig:result}) while the boiler is active for 11 hours in the steady-state scenario.
The reason is that for the given prices, the CHP is less expensive compared to the boiler.
Thus, the main cost reduction is achieved by shifting waste heat in time such that the boiler can be turned off.
This point is demonstrated by repeating the DR optimization but fixing the boilers on/off-decisions to the 11 hours of operation from the steady-state scenario.
The cost reduction is only 95.3 \euro, i.e., instead of 40.8\% improvement through DR only 5.7\% is realized.

In our case study, the simultaneous scheduling of processes and multi-energy system provides significant cost reductions via DR.
To realize the possible cost reductions, the discrete on/off-decisions in the multi-energy system need to be considered during optimization.
Using our high-order dynamic ramping constraints, we achieve optimization runtimes that allow an online application while the economic result is only slightly worse compared to MINLP optimization.

\subsection{Comparison with first-order static ramping constraints}
\begin{figure}[ht]    
\centering
\footnotesize
\includegraphics{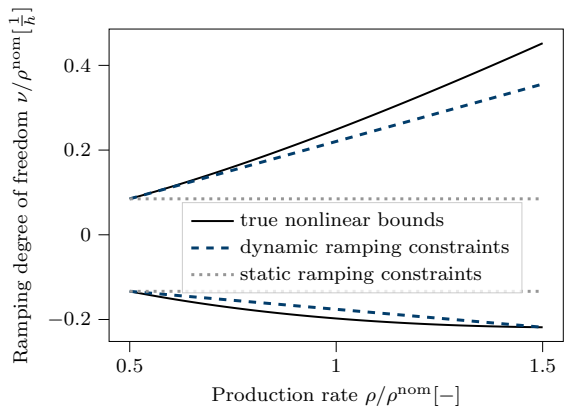}
    \caption{True nonlinear bounds, linear dynamic ramping constraints, and static ramping constraints for the non-jacketed CSTR 1 in the modified case study.}
    \label{fig:bounds50}
\end{figure}

Finally, we compare our high-order dynamic ramping constraints to first-order static ramping constraints.
As discussed, first-order ramping constraints are not applicable to the jacketed CSTR2 and we thus focus on the non-jacketed CSTR 1.
For CSTR1, the case study result does not change significantly if we use static instead of dynamic ramping constraints which is for two reasons: First, the difference between static and dynamic ramping constraints is small as we only vary the operating point by +/-~20\% (Figure \ref{fig:bounds}). 
Second, the optimal schedule shifts just enough waste heat in each time step such that the boiler can be turned off and thus the full ramping capabilities of CSTR 1 are not exploited.
To demonstrate the benefits of dynamic ramping constraints more clearly, we expand the operating range to +/- 50 \% such that the difference between dynamic and static ramping constraints becomes more pronounced (Figure \ref{fig:bounds50}).
Further, to study the ramping of CSTR 1 decoupled from on/off decisions, 
we modify the case study setup such that it contains only CSTR1, a CHP, and a constant heat demand. 
All produced electricity is sold to the grid. 
Moreover, we use a price profile that exhibits low prices for several hours and thus motivates a long period of ramping down. 
The price profile occurred on the 24th February 2018 at the German day-ahead market \citep{SMARD}.

\begin{figure}[h!]    
\centering
\includegraphics{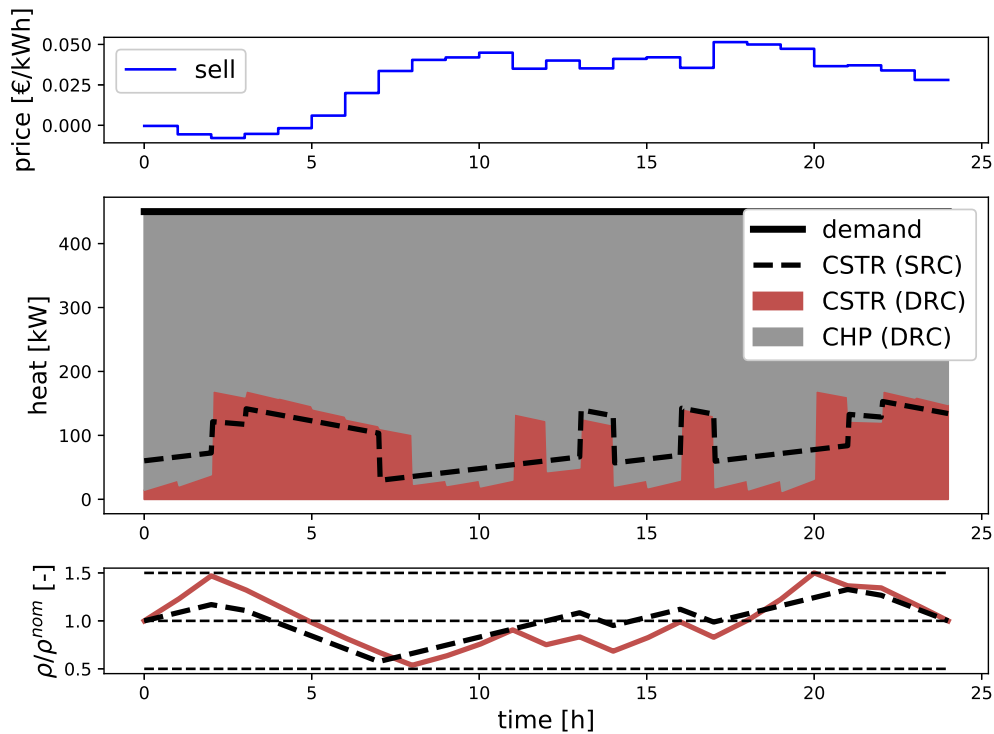}
    \caption{Resulting operation in the modified case study showing the heat supplied by CSTR and CHP, and the production $\rho$ delivered by the CSTR for the case of dynamic ramping constraints (DRC) (red line). For comparison, the waste heat and production rate with static ramping constraints (SRC) are shown as black dashed lines.}
    \label{fig:result_mod}
\end{figure}

The DR improvement is 6.7\% with static ramping constraints (SRC) while it becomes 12.2 \% with dynamic ramping constraints (DRC) (Table \ref{tab:results_comp}), i.e., the dynamic ramping constraints nearly double the benefits of demand response.
For static ramping constraints, the highest realized production rate is 133\% of the nominal production rate and the minimal production rate is 58\%.
In contrast, the range is 150\%-54\% with dynamic ramping constraints.
Especially, in the first two hours, the DRC scenario ramps up to 150\% while the SRC scenario only reaches 117\% (Figure \ref{fig:result_mod}).
This faster ramping increases the DR improvement because the morning hours have cheap electricity prices, which make operation of the CHP unfavorable, and thus favor a high waste heat production.
Most  waste heat is produced when the reactor ramps down as the temperature needs to be lowered during down-ramping.
While SRC ramp down starts from 117\% production rate in hour 2, DRC ramp down starts from 150\% production and therefore generates more waste heat.
The modified case study demonstrates that in cases where fast ramping is required dynamic ramping constraints can significantly improve demand response compared to static ramping constraints.

\begin{table}[h]
    \caption{Cost reduction achieved by waste heat integration in modified case study.}
    \label{tab:results_comp}
    \begin{tabular}{lcc}
    Case & Cost reduction [\euro] & DR improvement [\%] \\
    \hline
    No DR  & 81.1 & -\\
    DR with SRC& 86.5 & 6.7\\
    DR with DRC& 91.0 & 12.2\\
    \hline  
    \end{tabular}
\end{table}

\section{Discussion of possible extensions}
\label{sec:discussion}

In this section, we discuss four possible extensions of the dynamic ramping method, its integration with control, and the applicability of the approach.  

Extension 1: In some cases, it might be favorable to vary the output $y$ with the production rate instead of holding it constant because varying the output gives additional flexibility for example to choose optimal steady-state operating points.
The output can be varied if a function $\pi$ is chosen that defines the output $y$ as an $n$ times differentiable function of the production rate $\rho$. 
Consequently, the desired values of output $y$ and its time derivatives are (compare to Equations (\ref{eq:y_0}) - (\ref{eq:y_n})):
\begin{align}
    y & \Should \pi(\rho) \\
    \dot{y} & \Should \frac{\partial \pi(\rho)}{\partial \rho} \dot{\rho} \\
    y^{(2)} & \Should \frac{\partial^2 \pi(\rho)}{\partial^2 \rho} \dot{\rho}^2 + \frac{\partial \pi(\rho)}{\partial \rho} \rho^{(2)} \\
    &\vdots \nonumber
\end{align}
Still, the Equation system (\ref{eq:system_gamma}) is defined and our dynamic ramping reformulation can in principle be applied.
Note that the equation system is more complex to solve if the output $y$ is not constant anymore. 

Extension 2: Bounds on states are not considered in our example, however, they can be considered in a straightforward manner, as all states $\vect{x}$ are given as a function $\boldsymbol{\Gamma}(\boldsymbol{\varphi})$ of the ramping state vector $\boldsymbol{\varphi}$ (compare to Section~\ref{subsec:derive_dyn_ramping}).
For example, for the jacketed CSTR 2, the jacket temperature $T_j$ is given as a function $\Gamma(\rho,\dot{\rho})$ of the production rate $\rho$ and its first derivative, $\dot{\rho}$.
If the equation $T_j=\Gamma(\rho,\dot{\rho})$ can be solved for the derivative $\dot{\rho}$ to $\dot{\rho}=\Gamma^{-1}(\rho,T_j)$, one can insert the bounds of the state $T_j$ into the function $\Gamma^{-1}$ and receive bounds on the derivative $\dot{\rho}$ as a function of the production rate $\rho$.
Thus, bounds on the state $T_j$ can be considered by adding a first-order ramping constraint ($\dot{\rho}^{\text{min}}(\rho ) \leq \dot{\rho} \leq \dot{\rho}^{\text{max}}(\dot{\rho} )$).

Extension 3: As the presented derivation of dynamic ramping constraints is restricted to SISO processes, we continued to work on dynamic ramping constraints. We were able to generalize the rigorous derivation to flat MIMO processes \citep{Baader.2022}. 
Note that for SISO processes, flatness is equivalent to exact input-state linearizability, which we use in the current article.
The derivation of dynamic ramping constraints for MIMO processes follows the same concept as for SISO processes and is similar to the current paper's derivation.

Extension 4: Here, we only derive dynamic ramping constraints rigorously for the case that the relative degree $r$ equals the number of states $n$ ($r=n$). 
The case $r>n$ is not relevant as for $r>n$, the output $y$ is not controllable with the input $u$. For the case $r<n$, internal dynamics occur \citep{Corriou.2018}. 
Consequently, the process state is not fully described by $y$ and its first $r-1$ time derivatives but there are $n-r$ internal states. 
Such processes are called input-output linearizable \citep{Corriou.2018}. 
For input-output linearizable processes, our method would still allow to determine the order of dynamic ramping constraints rigorously. 
However, the limits $\nu^{\text{min}}$, $\nu^{\text{max}}$ varying with the operation point $\boldsymbol{\varphi}$ would additionally depend on the unknown internal states. 
Possibly, the internal states could be treated as bounded uncertainties and robust limits that hold for all possible values of the internal states could be calculated \citep{BenTal.2009}.
Alternatively, internal states might be estimated based on the trajectory of $\boldsymbol{\varphi}$, similar to the estimation of unmeasured intrinsic states performed by  \cite{Lovelett.2020}.

Integration with control: Our scheduling with dynamic ramping constraints can be combined with any regulatory control.
The result of the scheduling is a trajectory for the production rate $\rho$.
Due to the derivation of dynamic ramping constraints, it is guaranteed that if the trajectory of the production rate is applied to the process, the regulatory control can choose a feasible input $u$ ( $u^{\text{min}} \leq u \leq  u^{\text{max}} $) that keeps the output $y$ at its nominal value (compare to Section~\ref{subsec:derive_dyn_ramping}). 
However, if the scheduling decision is to ramp the production rate of the process as fast as possible, the input $u$ will be at its bound.
If in such a situation the controlled output $y$ deviates from its nominal value due to disturbances, noise, or model-plant mismatch, the control might be unable react as the input $u$ is already at its bound. 
Thus, a deviation of the output $y$ from its nominal value must be tolerated.
To solve this problem, the derivations of dynamic ramping constraints can be adapted to maintain the input within $u^{\text{min}} - \Delta u^{\text{control}} \leq u \leq  u^{\text{max}} - \Delta u^{\text{control}}$, with a safety margin $\Delta u^{\text{control}}$ that the regulatory control can use to compensate for disturbances.  

Even though any regulatory control can be used in principle, a feedforward linearization control structure \citep{Hagenmeyer2003,Hagenmeyer2003b} is particularly promising as the feedforward part simplifies the task for the regulatory control.
That is, the nominal control input $u^\text{nom}$ needed to reject the disturbances introduced by varying the production rate is calculated using Equation~(4.n+1) and directly applied to the process such that the underlying control does not have to compensate for the disturbance by the production rate anymore.
To reject other disturbances that might act on the process, a simple tracking controller, e.g., a PID-controller, is added to stabilize the feedforward control  \citep{Hagenmeyer2003,Hagenmeyer2003b}.

Overall, the applicability of our method is strongly limited due to the strict assumptions.
Still, as our results show that the dynamic ramping method could bridge the gap between nonlinear process models and simplified process representations for real-time scheduling, further research regarding the discussed extensions seems promising.
Additionally, we argue that the number of scheduling-relevant dynamics is typically small \citep{Du.2015,Baldea.2014}.
Consequently, even if the mechanistic process model under consideration might not fulfill our assumptions, it might be possible to consider a reduced-order model for the slow scheduling relevant dynamics to which dynamic ramping constraints can be applied.
For example, detailed electrolyzer models are typical MIMO models with several dynamics \citep{Hoffmann.2021}. 
However, typically, only the slow temperature dynamic has to be considered on the hourly time scale relevant for demand response \citep{Simkoff.2020,BenjaminFlamm.2021}.
In a recent conference publication, we show that dynamic ramping constraints can be transferred to electrolyzers with slow temperature dynamics and increase demand response potential compared to quasi-steady-state scheduling \citep{Baader.2021b}.

\section{Conclusion}
\label{sec:conclusion}
In this paper, we propose high-order dynamic ramping constraints for the simultaneous demand response (DR) optimization of processes and their multi-energy systems.
These dynamic ramping constraints can be of high order and the ramping limits depend on the process state.
Process-state-dependent limits enable faster transitions than typical static ramping constraints. 
Based on the notion of exact linearization, we derive dynamic ramping constraints rigorously for the case of exact input-state linearizable single-input single-output (SISO) processes. 

In a case study, we consider two continuous stirred-tank reactors (CSTRs) with waste heat integration that are scheduled simultaneously with a multi-energy system. 
Deriving dynamic ramping constraints from the two CSTR models, we formulate an MILP optimization problem and improve the economic value of the waste heat by 41 \% compared to steady-state operation.
This benefit is close to the bound of 43\% obtained from nonlinear mixed-integer dynamic optimization (MIDO) with the original process model.
Importantly, the MILP formulation based on dynamic ramping constraints allows to solve the simultaneous DR optimization within seconds.

In a modified case study, we find that the DR improvement with dynamic ramping constraints is significantly higher than that with static ramping as the state-dependent ramping limits allow faster ramping.
Consequently, the proposed high-order dynamic ramping constraints allow to capture the dynamic flexibility of processes better than traditional ramping constraints and achieve optimization runtimes sufficiently fast for online optimization.
  
  \section*{Author contributions}
\textbf{Florian J. Baader}: Conceptualization, Methodology, Software, Investigation, Validation, Visualization, Writing - original draft. 
\textbf{Philipp Althaus}: Conceptualization, Writing – review \& editing. 
\textbf{André Bardow}: Funding acquisition, Conceptualization, Supervision, Writing – review \& editing. 
\textbf{Manuel Dahmen}: Conceptualization, Supervision, Writing – review \& editing.

  \section*{Declaration of Competing Interest}
    We have no conflict of interest.

    \section*{Acknowledgements}
    This work was supported by the Helmholtz Association under the Joint Initiative “Energy System Integration”. AB and FB also received support from the Swiss Federal Office of Energy through the project "SWEET PATHFNDR".

    \newpage
\twocolumn
\section*{Nomenclature} \label{sec:nomenclature}

\noindent\textbf{Abbreviations} \\
\noindent
\begin{tabularx}{\columnwidth}{lX}
CSTR & continuous stirred tank reactors\\
DR & demand response \\
DRC & dynamic ramping constraint \\
MIDO & mixed-integer dynamic optimization \\
MILP & mixed-integer linear program \\
MIMO & multi-input multi-output \\
MINLP & mixed-integer nonlinear program \\
NLP & nonlinear program \\
SISO & single-input single-output \\
SRC & static ramping constraint
\end{tabularx} \\

\noindent\textbf{Greek symbols} \\
\noindent
\begin{tabularx}{\columnwidth}{lX}
$\alpha$ & nonlinear function \\
$\alpha_c$ & heat transfer coefficient \\
$\beta$ & nonlinear function \\
$\Gamma$ & state transformation\\
$\delta$ & order of DRC \\
$\eta$ & efficiency \\
$\nu$ & ramping degree of freedom  \\
$\rho$ & production rate \\
$\Phi$ & objective \\
$\tau$ & time constant \\
$\boldsymbol{\varphi}$ & ramping state vector \\
$\chi$ & optimization variable 
\end{tabularx} \\

\newpage

\noindent\textbf{Latin symbols} \\
\noindent
\begin{tabularx}{\columnwidth}{lX}
$a$ & fitting coefficient \\ 
$c$ &concentration \\
$F$ & flow rate \\
$f$ & nonlinear function \\
$g$ &nonlinear function \\
$h$ &nonlinear function \\
$\vect{J}$ & Jacobian matrix \\
$k$ & preexponential factor \\
$N$ & scaled activation energy \\
$n$ & number of states \\
$p$ &price \\
$Q$ & heat flow \\
$r$ & relative degree \\
$S$ & storage filling level \\
$T$ & temperature \\
$t$ & time \\
$u$ & input \\
$V$ & volume \\
\end{tabularx}
\begin{tabularx}{\columnwidth}{lX}
$x$ & differential state \\
$y$ & output \\
$z$ & discrete variable
\end{tabularx}

\noindent\textbf{Sets} \\
\noindent
\begin{tabularx}{\columnwidth}{lX}
$\mathbb{C}$ & energy system components \\
$\mathbb{C}_e^{\text{cons}}$ & energy system components that consume $e$\\
$\mathbb{C}_e^{\text{sup}}$ & energy system components that supply $e$ \\
$\mathbb{E}$ & end energy forms 
\end{tabularx} \\

\noindent\textbf{Subscripts} \\
\noindent
\begin{tabularx}{\columnwidth}{lX}
0 & initial value  \\
$c$ & cooling \\
dem & demand \\
$e$ & energy form \\
$f$ & feed \\
$i$ & energy system component \\
$j$ &jacket \\
$m$ &linear slope \\
$src$ & static ramping constraint \\
$wh$ & waste heat \\
\end{tabularx} \\

\noindent\textbf{Superscripts} \\
\noindent
\begin{tabularx}{\columnwidth}{lX}
$l$ & lower bound \\
max &maximum value \\
min &minimum value \\
nom & nominal value \\
$u$ & upper bound
\end{tabularx} \\

\newpage 
\onecolumn

  \bibliographystyle{apalike}
  \renewcommand{\refname}{Bibliography}  
  \bibliography{literature.bib}

\end{document}


\thispagestyle{firststyle}

  \begin{center}
    \begin{large}
      \textbf{\mytitle}
    \end{large} \\
    \myauthor
  \end{center}

  \vspace{0.5cm}

  \begin{footnotesize}
    \affil
  \end{footnotesize}

  \vspace{0.5cm}

\section*{Supplement: Complete derivatives of jacketed CSTR 2}
\label{ap:derivatives_long}
In equations (\ref{eq:CSTR2_h_A}) - (\ref{eq:CSTR2_y3_A}), we show the complete derivatives for the jacketed CSTR 2 (compare to Equations (14) - (17) in the main paper).

    \begin{align}
        \label{eq:CSTR2_h_A}
        y&=c \coloneqq \alpha_0(c) \Should c^{\text{nom}} \\
        \dot{y} &= \frac{\partial \alpha_0(c)}{\partial c} \dot{c} = (1 - c)\frac{\rho}{V}- c k e^{- \frac{N}{T}}\coloneqq \alpha_1(c,T,\rho) \Should 0  \\
        y^{(2)} &= \frac{\partial \alpha_1(c,T,\rho)}{\partial c} \dot{c} + \frac{\partial \alpha_1(c,T,\rho)}{\partial T} \dot{T} + \frac{\partial \alpha_1(c,T,\rho)}{\partial \rho} \dot{\rho} 
        \\ \nonumber
        &= - \left[\frac{\rho}{V} + k e^{- \frac{N}{T}}\right]\left[(1 - c)\frac{\rho}{V}- c k e^{- \frac{N}{T}}\right] 
        \\ \nonumber
        & - \left[  \frac{c k N e^{- \frac{N}{T}}}{T^2}\right]\left[(T_{f} - T)\frac{\rho}{V} + c k e^{- \frac{N}{T}} + \tau_1(T_j - T) \right] + \left[\frac{1-c}{V} \right] \dot{\rho} 
        \\ & \nonumber\coloneqq \alpha_2(c,T,T_j,\rho,\dot{\rho}) \Should 0 \\
        \label{eq:CSTR2_y3_A}
        y^{(3)} &= \frac{\partial \alpha_2(c,T,T_j,\rho,\dot{\rho})}{\partial \vect{x}} \left(\begin{array}{c} \dot{c} \\ \dot{T} \\ \dot{T}_j \end{array} \right) + \frac{\partial \alpha_2(c,T,T_j,\rho,\dot{\rho})}{\partial \boldsymbol{\varphi}} \left(\begin{array}{c} \dot{\rho} \\ \rho^{(2)}  \end{array} \right) 
        \\ \nonumber
        & = \frac{\partial \alpha_2(c,T,T_j,\rho,\dot{\rho})}{\partial c}  \left[(1 - c)\frac{\rho}{V}- c k e^{- \frac{N}{T}}\right]
        \\ \nonumber &+ \frac{\partial \alpha_2(c,T,T_j,\rho,\dot{\rho})}{\partial T}\left[(T_{f} - T)\frac{\rho}{V} + c k e^{- \frac{N}{T}} + \tau_1(T_j - T) \right]
        \\ \nonumber & -  \left[ \frac{ N \tau_1 c k  e^{- \frac{ N}{T}}}{T^{2}} \right]  \left[\tau_2(T - T_j) - F_{c} \alpha_c \left(T_j - T_{c} \right)\right]
        \\ \nonumber &+  \left[- \frac{ N c k \left(T_{f} - T\right) e^{- \frac{ N}{T}}}{T^{2} V} + \frac{\left(1 - c\right) \left(- k e^{- \frac{ N}{T}} - \frac{\rho}{V}\right)}{V} - \frac{- c k e^{- \frac{ N}{T}} + \frac{\rho \left(1 - c\right)}{V}}{V}\right] \dot{\rho}
        \\ \nonumber & + \left[\frac{ 1 - c}{V} \right] \rho^{(2)}
        \\ \nonumber
        &= \alpha_3(\vect{x},\boldsymbol{\varphi}) + \beta_u (\vect{x},\boldsymbol{\varphi}) F_c + \beta_{\rho}(\vect{x},\boldsymbol{\varphi})\nu \Should 0
    \end{align}
    \begin{align}
         \nonumber & \text{with } \nu = \rho^{(2)}, \frac{\partial \alpha_2(c,T,T_j,\rho,\dot{\rho})}{\partial c} =  - \frac{ N c k^{2} e^{- \frac{2  N}{T}}}{T^{2}} 
        \\ \nonumber 
        &- \frac{ N k \left( c k e^{- \frac{ N}{T}} - \tau_1 \left(- T_{j } + T\right) + \frac{\rho \left(T_{f} - T\right)}{V}\right) e^{- \frac{ N}{T}}}{T^{2}} + \left(- k e^{- \frac{ N}{T}} - \frac{\rho}{V}\right)^{2} - \frac{\dot{\rho}}{V}
        \\ \nonumber &
        \frac{\partial \alpha_2(c,T,T_j,\rho,\dot{\rho})}{\partial T} = - \frac{ N^{2} c k \left( c k e^{- \frac{ N}{T}} - \tau_1 \left(- T_{j } + T\right) + \frac{\rho \left(T_{f} - T\right)}{V}\right) e^{- \frac{ N}{T}}}{T^{4}}
        \\\nonumber 
        &- \frac{ N c k \left(- k e^{- \frac{ N}{T}} - \frac{\rho}{V}\right) e^{- \frac{ N}{T}}}{T^{2}}  - \frac{ N c k \left(\frac{ N c k e^{- \frac{ N}{T}}}{T^{2}} - \tau_1 - \frac{\rho}{V}\right) e^{- \frac{ N}{T}}}{T^{2}} 
        \\ \nonumber &- \frac{ N k \left(- c k e^{- \frac{ N}{T}} + \frac{\rho \left(1 - c\right)}{V}\right) e^{- \frac{ N}{T}}}{T^{2}} + \frac{2  N c k \left( c k e^{- \frac{ N}{T}} - \tau_1 \left(- T_{j } + T\right) + \frac{\rho \left(T_{f} - T\right)}{V}\right) e^{- \frac{ N}{T}}}{T^{3}} 
    \end{align}